\documentclass[12pt,a4paper]{amsart}
\usepackage{amssymb}
\usepackage{xypic}
\usepackage{rotating}
\newtheorem{thm}{Theorem}[section]

\newtheorem{prop}[thm]{Proposition}
\theoremstyle{definition}
\newtheorem{defin}[thm]{Definition}

\theoremstyle{remark}
\newtheorem{remark}[thm]{Remark}
\newtheorem{example}[thm]{Example}

\numberwithin{equation}{section}

\newcommand\Img{\operatorname{\mathcal Im}}
\newcommand\Hom{\operatorname{Hom}}

\newcommand\iso{\kern.35em{\raise3pt\hbox{$\sim$}\kern-1.1em\to}\kern.3em}

\newcommand\F{\mathbb F}
\newcommand\A{\mathbb A}
\newcommand\Pp{\mathbb P}
\newcommand\m{\mathfrak m}
\newcommand\Cc{\mathcal C}

\newcommand\Oc{\mathcal O}
\newcommand\Xcc{\mathcal X}

\newcommand\Res{\operatorname {Res}}
\newcommand\Spec{\operatorname {Spec}}
\newcommand\Proj{\operatorname {Proj}}

\begin{document}
\title{Convolutional Goppa Codes}
\author{J.M. Mu\~noz Porras, {\sl Member, IEEE}}
\author{J.A. Dom\'{\i}nguez P\'erez}
\author{J.I. Iglesias Curto}
\author{G. Serrano Sotelo}
\email{jmp@usal.es, jadoming@usal.es, joseig@usal.es and
laina@usal.es}
\address{Departamento de Matem\'aticas, Universidad de Salamanca,
Plaza de la Merced 1-4, 37008 Salamanca, Spain}
\date{\today}
\thanks {This research was partially supported by the Spanish DGESYC
through research project BMF2000-1327 and by the ``Junta de
Castilla y Le\'on'' through research projects SA009/01 and
SA032/02.}
%\keywords{}

\begin{abstract}
We define Convolutional Goppa Codes over algebraic curves and
construct their corresponding dual codes. Examples over the
projective line and over elliptic curves are described, obtaining
in particular some Maximum-Distance Separable (MDS) convolutional
codes.
\end{abstract}

\maketitle
\begin{quote}
\end{quote}

\section{Introduction}

Goppa codes are evaluation codes for linear series over smooth
curves over a finite field $\F_q$. In \cite{DMS} we proposed a new
construction of convolutional codes, which we called Convolutional
Goppa Codes (CGC), in terms of evaluation along sections of a
family of algebraic curves.

The aim of this paper is to reformulate the results of \cite{DMS}
in a straightforward language. We define CGC as Goppa codes for
smooth curves defined over the field $\F_q(z)$ of rational
functions in one variable $z$ over the finite field $\F_q$. These
CGC are in fact more general than the codes defined in \cite{DMS},
since there are smooth curves over $\F_q(z)$ that do not extend to
a family of smooth curves over the affine line $\A_{\F_q}^1$. With
this definition, one has another advantage: the techniques of
Algebraic Geometry we need are easier than those used in
\cite{DMS}: we use exactly the same language as is usual in the
literature on Goppa codes.

The last two sections of the paper are devoted to illustrating the
general construction with some examples. In \S 4 we construct
several CGC of genus zero; that is, defined in terms of the
projective line $\Pp_k^1$ over the field $\F_q(z)$. Some of these
examples are MDS-convolutional codes and are very easy to handle.

In \S 5 we give examples of CGC of genus one; that is, defined in
terms of elliptic curves over $\F_q(z)$. These examples are not so
easy to study. In fact, a consequence of this preliminary study of
CGC of genus one is that a deeper understanding of the arithmetic
properties of elliptic fibrations (see for instance \cite{Ta}) and
of the translation of these properties into the language of
convolutional codes, is necessary.

\section{Convolutional Goppa Codes}

Let $\F_q$ be a finite field and $\F_q(z)$ the (infinite) field of
rational functions of one variable. Let $(X,\Oc_X)$ be a smooth
projective curve over $\F_q(z)$ of genus $g$, and let us denote by
$\Sigma_X$ the field of rational functions of $X$.

Given a set $p_1,\dots, p_n$ of $n$ different $\F_q(z)$-rational
points of $X$,  if $\Oc_{p_i}$ denotes the local ring at the point
$p_i$, with maximal ideal $\m_{p_i}$, and $t_i$ a local parameter
at $p_i$, one has exact sequences
\begin{equation}\label{point}
\begin{aligned}
0\to \m_{p_i} \to \Oc_{p_i}&\to \Oc_{p_i}/\m_{p_i}\simeq
\F_q(z)\to 0\cr s(t_i)&\mapsto s(p_i)
\end{aligned}
\end{equation}
Let us consider the divisor $D=p_1+\dots+p_n$, with its associated
invertible sheaf $\Oc_X(D)$. Then, one has an exact sequence of
sheaves
\begin{equation}\label{divisor}
0\to \Oc_X(-D)\to \Oc_X\to Q\to 0\,,
\end{equation}
where the quotient $Q$ is a sheaf with support at the points
$p_i$.

Let $G$ be a divisor on $X$ of degree $r$, with support disjoint
from $D$. Tensoring the exact sequence (\ref{divisor}) by the
associated invertible sheaf $\Oc_X(G)$, one obtains:
\begin{equation}\label{sheaf}
0\to \Oc_X(G-D)\to \Oc_X(G)\to Q\to 0\,.
\end{equation}

For every divisor $F$ over $X$, let us denote their
$\F_q(z)$-vector space of global sections by
$$L(F)\equiv
\Gamma(X,\Oc_X(F))=\{s\in \Sigma_X\ /\ (s)+F\geq 0\}\,,$$ where
$(s)$ is the divisor defined by $s\in \Sigma_X$. Taking global
sections in (\ref{sheaf}), one obtains
$$
\begin{aligned}
0\to L(G-D) \to L(G)& \overset\alpha\to
\F_q(z)\times\overset{n)}\dots\times \F_q(z)\to \dots\cr s&\mapsto
(s(p_1),\dots, s(p_n))
\end{aligned}
$$

\begin{defin}\label{cgc} The convolutional Goppa code $\Cc(D,G)$ associated
with the pair $(D,G)$ is the image of the $\F_q(z)$-linear map
$\alpha\colon L(G)\to \F_q(z)^n$.

Analogously, given a subspace $\Gamma\subseteq L(G)$, one defines
the convolutional Goppa code $\Cc(D,\Gamma)$ as the image of
${\alpha}_{\vert_{\Gamma}}$.
\end{defin}

\begin{remark} The above definition is more general than the one given in
\cite{DMS} in terms of families of curves $X\to\A^1_{\F_q}$. In
fact, given such a family, the fibre $X_\eta$, over the generic
point $\eta \in \A^1_{\F_q}$, is a curve over $\F_q(z)$. But not
every curve over $\F_q(z)$ extends to a family over $\A^1_{\F_q}$.
\end{remark}

By construction, $\Cc(D,G)$ is a convolutional code of length $n$
and dimension
$$k\equiv\dim L(G)-\dim L(G-D)\,.$$

\begin{prop}\label{cgc} Let us assume that $2g-2<r<n$. Then, the
evaluation map $\alpha\colon L(G)\hookrightarrow \F_q(z)^n$ is
injective, and the dimension of $\Cc(D,G)$ is
$$k=r+1-g\,.$$
\end{prop}
\begin{proof} If $r<n$, $\dim L(G-D)=0$, the map $\alpha$ is injective and
$k=\dim L(G)$. If $2g-2<r$, $\dim L(G)=1-g+r$ by the Riemann-Roch
theorem.
\end{proof}

\section{Dual Convolutional Goppa Codes}

Let us consider, over the $\F_q(z)$-vectorial space $\F_q(z)^n$,
the pairing $\langle\ ,\ \rangle$
$$
\begin{aligned}
\F_q(z)^n \times &\F_q(z)^n  \to \F_q(z)\cr (u,&v)\longmapsto
\langle u , v \rangle = \sum_{i=1}^n u_i v_i\,,
\end{aligned}
$$
where $u=(u_1,\dots,u_n), v=(v_1,\dots,v_n)\in \F_q(z)^n$.

 \begin{defin}\label{cgcd}
The dual convolutional Goppa code of the code $\Cc(D,G)$ is the
$\F_q(z)$-linear subspace $\Cc^\perp(D,G)$ of $\F_q(z)^n$ given by
$$\Cc(D,G)^\perp=\{ u\in \F_q(z)^n \ / \ \langle u , v \rangle = 0 \text{\ for every\
} v\in \Cc(D,G)\}\,.
$$
\end{defin}

Let us denote by $K$ the canonical divisor of rational
differential forms over $X$.

\begin{thm}\label{dcgc} The dual convolutional Goppa code $\Cc^\perp(D,G)$ associated
with the pair $(D,G)$ is the image of the $\F_q(z)$-linear map
$\beta\colon L(K+D-G)\to \F_q(z)^n$, given by
$$\beta(\eta)=(\Res_{p_1}(\eta),\dots, \Res_{p_n}(\eta))\,.$$
\end{thm}

\begin{proof}
Following the construction of $\Cc(D,G)$, we start tensoring the
exact sequence (\ref{point}) by
$\m_{p_i}^\ast=\Hom_{\Oc_{p_i}}(\m_{p_i},\Oc_{p_i})$, and obtain
\begin{equation}\label{dpoint}
\begin{aligned}
0\to \Oc_{p_i} \to \m_{p_i}^\ast&\longrightarrow
\Oc_{p_i}/\m_{p_i}\otimes_{\Oc_{p_i}}\m_{p_i}^\ast\simeq
\F_q(z)\to 0\cr t_i^{-1}s(t_i&)\mapsto s(p_i)
\end{aligned}
\end{equation}
Again tensoring (\ref{dpoint}) by $\m_{p_i}/\m_{p_i}^2$, the
tangent space of differentials at the point $p_i$, one obtains
\begin{equation}\label{ddpoint}
\begin{aligned}
0\to \m_{p_i}/\m_{p_i}^2 \to
\m_{p_i}^\ast\otimes_{\Oc_{p_i}}\m_{p_i}/&\m_{p_i}^2\longrightarrow
\F_q(z)\longrightarrow 0\cr t_i^{-1}s(t_i)dt_i&\mapsto
s(p_i)=\Res_{p_i}(t_i^{-1}s(t_i)dt_i)\,.
\end{aligned}
\end{equation}

This allows us to define a new convolutional Goppa code associated
to the pair of divisors $D=p_1+\dots+p_n$ and $G$; tensoring
(\ref{divisor}) by the line sheaf $\Oc_X(K+D-G)$, one has
\begin{equation}\label{sheaf}
0\to \Oc_X(K-G)\to \Oc_X(K+D-G)\to Q\to 0\,.
\end{equation}
Taking global sections, one has:
$$
\begin{aligned}
0\to L(K-G) \to L(K+D-G)& \overset\beta\to
\F_q(z)\times\overset{n)}\dots\times \F_q(z)\to \dots\cr
\eta&\mapsto (\Res_{p_1}(\eta),\dots, \Res_{p_n}(\eta))
\end{aligned}
$$
The image of $\beta$ is a subspace of $\F_q(z)^n$, whose dimension
can be calculated by the Riemann-Roch theorem:
$$\begin{aligned}\dim L(K+&D-G)-\dim L(K-G)=\cr=&(\dim L(G-D)-(r-n)-1+g)-(\dim
L(G)-r-1+g)=n-k\,.
\end{aligned}
$$
Moreover, $\Img \beta$ is the subspace $\Cc(D,G)^\perp\subset
\F_q(z)^n$, since they have the same dimension, and for every
$\eta\in L(K+D-G)$ and every $s\in L(G)$ one has
$$\langle \beta(\eta), \alpha(s)\rangle =\sum_{i=1}^n s(p_i)\Res_{p_i}(\eta)=\sum_{i=1}^n
\Res_{p_i}(s\eta)=0\,,
$$
by the Residue Theorem.
\end{proof}

Under the hypothesis $2g-2<r<n$, the map $\beta$ is injective, and
$\Cc^\perp(D,G)$ is a convolutional code of length $n$ and
dimension
$$\dim L(K+D-G)=n-(1-g+r)\,.$$

\begin{remark} Our pairing $\langle\ ,\ \rangle\colon \F_q(z)^n \times \F_q(z)^n  \to
\F_q(z)$ is $\F_q(z)$-bilinear, whereas the ``time reversal''
pairing defined by J. Rosenthal in \cite{Ro} $7.5$, given by
$$
\begin{aligned}
\left[\ ,\ \right]\colon \F_q((z))^n \times &\F_q((z))^n  \to
\F_q\cr (u,&v)\longmapsto \sum_{i=1}^n \langle u(i) ,
v(-i)\rangle\,,
\end{aligned}
$$
where $u=\sum_i u(i)z^i, v=\sum_i v(i)z^i\in \F_q((z))^n$ and
$\langle\ ,\ \rangle$ is the standard bilinear form on $\F_q^n$,
 is $\F_q$-bilinear.

 The pairing $\left[\ ,\ \right]$ can be expressed in the following way:
$$\left[ u, v\right] =\Res_{z=0}\left(\langle u,v\rangle
\frac{dz}{z}\right)=\sum_{i=1}^n \Res_{z=0}\left(u_i v_i
\frac{dz}{z}\right) \,.$$ Thus, the duality for convolutional
Goppa codes defined in {\bf \ref{cgcd}} is related to the residues
in the points of $X$, and the duality with respect to the pairing
$\left[\ ,\ \right]$ is related to the residues in the variable of
the base field. A more precise study of the relationship between
both dualities must be done.
\end{remark}

\section{Convolutional Goppa Codes over the projective line}

Let $X=\Pp^1_{\F_q(z)}=\Proj \F_q(z)[x_0,x_1]$ be the projective
line over the field $\F_q(z)$, and let us denote by $t=x_1/x_0$
the affine coordinate.

Let $p_0=(1,0)$ be the origin point, $p_\infty=(0,1)$ the point at
infinity, and $p_1,\dots,p_n$ be different rational points of
$\Pp^1$, $p_i\neq p_0, p_\infty$. Let us define the divisors
$D=p_1+\dots+p_n$ and $G=r p_\infty - s p_0 $, with
$$0\leq s \leq r<n\,.$$
Since $g=0$, the evaluation map $\alpha\colon L(G)\to \F_q(z)^n$
is injective and $\Img \alpha$ defines a convolutional Goppa code
$\Cc(D,G)$ of length $n$ and dimension $k=r-s+1$.

Let us choose the functions $t^s,t^{s+1},\dots,t^r$ as a basis of
$L(G)$. If $\alpha_i\in\F_q(z)$ is the local coordinate of the
point $p_i$, $i=1,\dots,n$, the matrix of the evaluation map
$\alpha$ is,
\begin{equation}\label{matrixg}
G=\left( \begin{matrix} \alpha_1^s & \alpha_2^s & \dots &
\alpha_n^s\cr \alpha_1^{s+1} & \alpha_2^{s+1} &\dots &
\alpha_n^{s+1} \cr \vdots & \vdots &\ddots & \vdots\cr \alpha_1^r
& \alpha_2^r & \dots & \alpha_n^r\cr
\end{matrix}\right)\,.
\end{equation}
This is a generator matrix for the code $\Cc(D,G)$.

The dual convolutional Goppa code $\Cc^\perp(D,G)$ also has length
$n$, and dimension $n-k=n-r+s-1$. To construct $\Cc^\perp(D,G)$,
let us choose in $L(K+D-G)$ the basis of rational differential
forms $\left\langle {\frac{dt}{t^s\prod_{i=1}^n(t-\alpha_i) }},
{\frac{t\ dt}{t^s\prod_{i=1}^n(t-\alpha_i) }}, \dots,
{\frac{t^{n-r+s-2} dt}{t^s\prod_{i=1}^n(t-\alpha_i)
}}\right\rangle$, and let us calculate the residues
$$\begin{aligned}
\Res_{p_j}\!\left(\frac{t^m dt}{t^s\prod_{i=1}^n(t-\alpha_i)
}\right)= &\Res_{p_j}\!\left( \frac{(t-\alpha_j+\alpha_j)^m
d(t-\alpha_j)}
{(t-\alpha_j)(t-\alpha_j+\alpha_j)^s\prod_{\overset{i=1}{i\neq
j}}^n(t-\alpha_j+\alpha_j-\alpha_i)} \right)=\cr
=&{\frac{\alpha_j^m}{\alpha_j^s\prod_{\overset{i=1}{i\neq
j}}^n(\alpha_j-\alpha_i) }}
\end{aligned}
$$
If one denotes by
$h_j=\frac{1}{\alpha_j^s\prod_{\overset{i=1}{i\neq
j}}^n(\alpha_j-\alpha_i) }$, then the matrix of $\beta\colon
L(K+D-G)\to \F_q(z)^n$
\begin{equation}\label{matrixc}
H=\left( \begin{matrix} h_1  & h_2 &\dots & h_n\cr h_1\alpha_1 &
h_2\alpha_2 &\dots & h_n\alpha_n\cr \vdots & \vdots &\ddots &
\vdots\cr h_1\alpha_1^{n-r+s-2} & h_2\alpha_2^{n-r+s-2} & \dots &
h_n\alpha_n^{n-r+s-2}\cr
\end{matrix}\right)\,,
\end{equation}
is a generator matrix for the dual code $\Cc^\perp(D,G)$, and
therefore a parity-check matrix for $\Cc(D,G)$. In fact, one has
$H\cdot G^T=0$.

\begin{remark} The matrix (\ref{matrixc}) suggests that
$\Cc^\perp(D,G)$ is an alternant code over the field $\F_q(z)$,
and we can thus apply some kind of Berlekamp-Massey decoding
algorithm for convolutional Goppa codes; this will be studied in a
forthcoming paper.
\end{remark}

\begin{example}\label{ex} Let $a,b\in \F_q$ be two different non-zero
elements, and
$$\alpha_i=a^{i-1} z +
b^{i-1}\,, i=1,\dots, n\,, \text{ with $n<q$}\,. $$ We present
some examples of convolutional Goppa codes with canonical
generator matrices \cite{McE}, whose {\em free distance d\/}
attains the generalized Singleton bound, i.e., they are MDS
convolutional codes \cite{RoSm}, and we include their encoding
equations as linear systems
$$\left.\begin{aligned}
z^{-1}s=& s A_{\delta\times\delta} + u B_{k\times \delta}\cr u G=&
s C_{\delta\times n} + u D_{k\times n}
\end{aligned}\right\}\,,
$$
where $\delta$ denotes the degree of the code (in the sense of
\cite{McE}.)

\bigskip

\begin{itemize}
\item Field $\F_3(z)$, $\F_3=\{0,1,2\}$:

\medskip
$\begin{aligned} &G= \left( \begin{matrix} z+1 & z+2
\end{matrix} \right)\,,\cr
&H=\left( \begin{matrix} \frac{1}{2(z+1)}& \frac{1}{z+2}
\end{matrix} \right)\,,\cr
&A=\left( \begin{matrix} 0 \end{matrix} \right) \,,\quad
B=\left(\begin{matrix} 1 \end{matrix} \right)\,,\quad C=\left(
\begin{matrix} 1 & 1 \end{matrix} \right) \,,\quad
D=\left(
\begin{matrix} 1 & 2 \end{matrix} \right) \,,\cr
&(n,k,\delta,d)=(2,1,1,4)\,.\cr &
\end{aligned}
$

\item Field $\F_4(z)$, $\F_4=\{0,1,\alpha,\alpha^2\}$ where
$\alpha^2+\alpha+1=0$:

\medskip
$\begin{aligned} &G=\left( \begin{matrix} 1 & 1 & 1\cr z+1 &
\alpha z+\alpha^2 & \alpha^2 z+\alpha
\end{matrix} \right)\,,\cr
&H=\left( \begin{matrix} \frac{1}{(\alpha^2 z+\alpha)(\alpha
z+\alpha^2)} & \frac{1}{(\alpha^2 z+\alpha)(z+1)} &
\frac{1}{(\alpha z+\alpha^2)(z+1)}
\end{matrix} \right)\,,\cr
& A=\left( \begin{matrix} 0 \end{matrix} \right)\,,\quad B=\left(
\begin{matrix} 0 \cr 1 \end{matrix} \right)\,,\quad C=
\left( \begin{matrix} 1 & \alpha & \alpha^2\end{matrix}
\right)\,,\quad D=\left(
\begin{matrix} 1 & 1 & 1\cr 1 & \alpha^2 & \alpha\end{matrix} \right)\,,\cr
&(n,k,\delta,d)=(3,2,1,3)\,.\cr &
\end{aligned}
$

\item Field $\F_4(z)$:

\medskip
$\begin{aligned} &G=\left( \begin{matrix} z+1 & z+\alpha &
z+\alpha^2
\end{matrix} \right)\,,\cr
&H=\left( \begin{matrix} \frac{1}{z+1} & \frac{\alpha}{z+\alpha} &
\frac{\alpha^2}{z+\alpha^2}\cr 1 & \alpha & \alpha^2
\end{matrix} \right)\,,\cr
& A=\left( \begin{matrix} 0 \end{matrix} \right)\,,\quad B=\left(
\begin{matrix} 1 \end{matrix} \right)\,,\quad C=
\left( \begin{matrix} 1 & 1 & 1\end{matrix} \right)\,,\quad
D=\left(
\begin{matrix} 1 & \alpha & \alpha^2\end{matrix} \right)\,,\cr
&(n,k,\delta,d)=(3,1,1,6)\,.\cr &
\end{aligned}
$

\item Field $\F_5(z)$, $\F_5=\{0,1,2,3,4\}$:

\medskip
$\begin{aligned} &G=\left( \begin{matrix} (z+1)^2 & (z+2)^2 &
(z+4)^2
\end{matrix} \right)\,,\cr
&H=\left( \begin{matrix} \frac{2}{(z+1)^2} & \frac{2}{(z+2)^2} &
\frac{1}{(z+4)^2}\cr \frac{2}{z+1} & \frac{2}{z+2} & \frac{1}{
z+4}
\end{matrix} \right)\,,\cr
& A=\left( \begin{matrix} 0 & 1\cr 0& 0\end{matrix}
\right)\,,\quad B=\left(
\begin{matrix} 1 & 0 \end{matrix} \right)\,,\quad C=
\left( \begin{matrix} 2 & 4 & 3\cr 1 & 1 & 1\end{matrix}
\right)\,,\quad D=\left(
\begin{matrix} 1 & 4 & 1\end{matrix} \right)\,,\cr
&(n,k,\delta,d)=(3,1,2,9)\,.\cr &
\end{aligned}
$

\item Field $\F_5(z)$:

\medskip
$\begin{aligned} &G=\left( \begin{matrix} z+1 & 2z+3 & 4z+4 & 3z+2
\cr (z+1)^2 & (2z+3)^2 & (4z+4)^2 & (3z+2)^2
\end{matrix} \right)\,,\cr
&H=\left( \begin{matrix} \frac{4}{(z+1)^2(z+2)(z+3)} &
\frac{4}{(z+2)(z+3)(z+4)^2} & \frac{4}{(z+1)^2(z+2)(z+3)} &
\frac{4}{(z+2)(z+3)(z+4)^2} \cr \frac{4}{(z+1)(z+2)(z+3)} &
\frac{3}{(z+2)(z+3)(z+4)} & \frac{1}{(z+1)(z+2)(z+3)} & \frac{2}{
(z+2)(z+3)(z+4)}
\end{matrix} \right)\,,\cr
& A=\left( \begin{smallmatrix} 0 & 0 & 0\cr 0 & 0 & 1 \cr 0 & 0 &
0\end{smallmatrix} \right)\,,\quad B=\left(
\begin{smallmatrix} 1 & 0 & 0 \cr 0 & 1 &
0 \end{smallmatrix} \right)\,,\quad C= \left( \begin{smallmatrix}
1 & 2 & 4 & 3 \cr 2 & 2 & 2 & 2 \cr 1 & 4 & 1 & 4\end{smallmatrix}
\right)\,,\quad D=\left(
\begin{smallmatrix} 1 & 3 & 4 & 2 \cr 1 & 4 & 1 & 4\end{smallmatrix} \right)\,,\cr
&(n,k,\delta,d)=(2,1,3,8)\,.\cr &
\end{aligned}
$
\end{itemize}

\end{example}

\section{Convolutional Goppa Codes associated with elliptic curves}

    We can obtain convolutional codes from elliptic curves in the same way.
Let $X\subset \Pp^2_{\F_q(z)}$ be a plane elliptic curve over
$\F_q(z)$, and let us denote by $(x,y)$ the affine coordinates in
$\Pp^2_{\F_q(z)}$. Let $p_\infty$ be the infinity point, and
$p_1,\dots,p_n$ rational points of $X$, with $p_i =
(x_i(z),y_i(z))$. Let us define $D=p_1+\dots + p_n$ and $G=r
p_\infty$.

    The ``canonical" basis of $L(G)$ is $\{1,x,y,\ldots,x^ay^b\}$, with $2a+3b=r$.
Thus, the evaluation map $\alpha\colon L(G)\to \F_q(z)^n$ is
    \begin{equation*}
        \alpha(x^iy^j) =
        (x^i_1(z)y^j_1(z),\ldots,x^i_n(z)y^j_n(z)).
    \end{equation*}
The image of a subspace $\Gamma\subseteq L(G)$ under the map
$\alpha$ provides a Goppa convolutional code.

    We present a couple of examples obtained from elliptic curves that,
although not MDS, have free distance approaching that bound.

\begin{example}

     We consider the curve over $\F_2(z)$
    \begin{equation*}
    y^2 + (1+z)xy +(z+z^2)y = x^3 + (z+z^2) x^2
    \end{equation*}
    and the points
    \begin{equation*}
    \begin{array}{l}
        p_1 = (z^2 + z,z^3+z^2)\\
        p_2 = (0, z^2 + z)\\
        p_3 = (z, z^2)
    \end{array}
    \end{equation*}
    Let $\Gamma\subset L(G)$ be the
    subspace generated by $\{1,x\}$. Accordingly, the valuation map
    $\alpha$ over $\Gamma$ is defined by the matrix
    \begin{equation*}
    \left (
        \begin{array}{ccc}
        1 & 1 & 1 \\
        z^2+z & 0 & z
        \end{array}
    \right ).
    \end{equation*}
    This code has free distance $d=2$. The maximum distance for its
    parameters is $3$.

\end{example}

\begin{example}

     Let us now consider the curve over $\F_2(z)$
    \begin{equation*}
    y^2 + (1+z+z^2)xy +(z^2+z^3)y = x^3 + (z^2+z^3) x^2
    \end{equation*}
    and the points
    \begin{equation*}
    \begin{array}{l}
          p_1 = (z^3 + z^2, 0)\\
        p_2 = (0, z^3 + z^2)\\
        p_3 = (z^3 + z^2, z^5 + z^3)\\
        p_4 = (z^2 + z,z^3 + z)\\
        p_5 = (z^2 + z, z^4 + z^2)
    \end{array}
    \end{equation*}
    Let $\Gamma\subset L(G)$ be the
    subspace generated by $\{1,x\}$. Then, the valuation map
    $\alpha$ over $\Gamma$ is defined by the matrix
    \begin{equation*}
    \left (
        \begin{array}{ccccc}
        1 & 1 & 1 & 1 & 1 \\
        z^3+z^2 & 0 & z^3+z^2 & z^2 + z & z^2 + z
        \end{array}
    \right ).
    \end{equation*}
    This code has free distance $d=4$. The maximum distance for its
    parameters is $5$.

\begin{remark} Every elliptic curve $X$ over $\F_q(z)$ can be
considered as the generic fibre of a fibration $\Xcc\to U=\Spec
\F_q[z]$, with some fibres singular curves of genus $1$. The
global structure of this fibration is related to the singular
fibres (see \cite{Ta}); the translation into the language of
coding theory of the arithmetic and geometric properties of the
fibration is the first step in the program of applying the general
construction to the effective construction of good convolutional
Goppa codes of genus $1$.
\end{remark}

\end{example}

\medskip
\noindent {\bf Acknowledgments.} We thank J. Prada Blanco for
helpful comments that served to stimulate our work.


\begin{thebibliography}{99}



\bibitem{DMS}
J.A. Dom\'{\i}nguez P\'erez, J.M.Mu\~noz Porras and G. Serrano
Sotelo, Convolutional Codes of Goppa Type, {\em submitted to
AAECC\/}.


\bibitem{McE}
R.J. McEliece, The Algebraic Theory of Convolutional Codes, in:
{\em Handbook of Coding theory\/}, Ed. by V.S. Pless and W.C.
Huffman (Elsevier, Amsterdam, 1998) 1065--1138.


\bibitem{Ro}
J. Rosenthal, Connections between linear systems and convolutional
codes, in: {\em Codes, Systems and Graphical Models\/}, Ed. by B.
Marcus and J. Rosenthal, IMA, vol. 123, (Springer-Verlag, 2000)
39--66.


\bibitem{RoSm}
J. Rosenthal and R. Smarandache, Maximum Distance Separable
Convolutional Codes, {\em AAECC\/} {\bfseries 10} (1999) 15--32.


\bibitem{Ta}
J. Tate, Algorithm for determining the type of a singular fiber in
an elliptic pencil, in: {\em Modular Functions of One Variable\/},
LNM, vol. 476, (Springer-Verlag, 1975).


\end{thebibliography}
\end{document}